\documentclass[a4paper,12pt,reqno]{amsart}
\usepackage[utf8]{inputenc}
\usepackage{amssymb}
\usepackage{amsmath}
\usepackage{makecell}
\usepackage{multirow}
\usepackage{hhline}
\usepackage{amsthm}
\usepackage{amsfonts}
\usepackage{mathtools}
\usepackage{todonotes}
\usepackage{comment}
\usepackage{mathrsfs}
\usepackage{enumitem}
\usepackage[all,cmtip]{xy}
\usepackage[hyphens]{url}
\usepackage{hyperref}
\usepackage{cleveref}

\usepackage[british]{babel}
\usepackage[margin=1.2in]{geometry}
\usepackage{eqparbox}

\usepackage{pbox}

\renewcommand{\O}{\mathcal{O}} 

\newcommand{\set}[1]{\left\lbrace #1 \right\rbrace}

\newcommand{\field}[1]{\mathbb{#1}}  
\newcommand{\Q}{\field{Q}} 
\newcommand{\C}{\field{C}} 

\renewcommand{\P}{\field{P}}
\newcommand{\PP}{\field{P}}

\DeclareMathOperator{\Contr}{Contr}

\DeclareMathOperator{\Pic}{Pic}

\DeclareMathOperator{\Gal}{Gal}

\newtheorem{lemma}{Lemma}
\newtheorem{theorem}[lemma]{Theorem}
\newtheorem{proposition}[lemma]{Proposition}
\newtheorem{corollary}[lemma]{Corollary}

\theoremstyle{definition}
\newtheorem{definition}[lemma]{Definition}

\newtheorem{remark}[lemma]{Remark}

\numberwithin{lemma}{section}
\numberwithin{equation}{section} 
\numberwithin{figure}{section}
\title{Large degree primitive points on curves}
\author{Maarten Derickx}

\begin{document}

\begin{abstract}
A number field $K$ is called primitive if $\Q$ and $K$ are the only subfields of $K$. Let $X$ be a nice$^1$\thanks{$^1$nice means smooth projective and geometrically irreducible.} curve over $\Q$ of genus $g$. A point $P$ of degree $d$ on $X$ is called primitive if the field of definition $\Q(P)$ of the point is a primitive extension of $\Q$. In this short note we prove that if $X$ has a divisor of degree $d> 2g$, then $X$ has infinitely many primitive points of degree $d$. This complements recent results of Khawaja and Siksek in \cite{KhawajaSiksek} that show that under certain conditions most of the points of low degree are not primitive.
\end{abstract}

\maketitle

\section{Introduction}

Khawaja and Siksek studied the set of primitive points of a curve. One of the main results on this is the following:

\begin{theorem}[{\cite{KhawajaSiksek}[Thm. 2]}] Let $X$ be a nice curve over $\Q$ with genus $g$ and $\Q$-gonality $m \geq 2$.
Let $d \geq 2$ be an integer satisfying \begin{align}d \neq m, \quad d < 1+ \frac {g} {m-1} \label{eq:degree_condition}.\end{align} Let $J$ be the Jacobian of $X$. Suppose either
of the following hold:
\begin{itemize}
    \item[(a)] $J(\Q)$ is finite;
    \item[(b)] or $d \leq g - 1$, and $A(\Q)$ is finite for every abelian subvariety $A/\Q$ of $J$ of dimension $\leq d/2$.
\end{itemize}
    
Then $X$ has finitely many primitive degree d points. Moreover, if $\gcd(d, m) = 1$ or
$d$ is prime then $X$ has finitely many degree $d$ points.
\end{theorem}

Let $d$ be an integer satisfying \cref{eq:degree_condition}. The argument used in\cite{KhawajaSiksek} relies on the fact that all functions $f$ on degree $d$ on $X$ have to factor in a nontrivial way as composition of two morphisms $X \to C \to \P^1$ due to the Castelnuovo–Severi inequality. However if $d$ is large enough with respect to $g$ then $X$ will always have an abundance of functions of degree $d$, as soon as it has at least one function of degree $d$. Additionally there is no theorem anymore forcing all these functions to factor trough some intermediate curve. In particular, the argument that forces all but finitely many points to be imprimitive breaks down for larger degree.  So it is reasonable to expect an abundance of primitive points for larger degrees. 

In the same article Khawaja and Siksek already give an example that shows that the finiteness of primitive points breaks down when the degree is large enough.

\begin{lemma}[{\cite{KhawajaSiksek}[Lem. 31]}]
Let $g \geq 2$ and  $d = g + 1$. Then there is a hyperelliptic curve $C/\Q$ of genus $g$ with infinitely many primitive degree d points.
\end{lemma}

The main result of this article is that the finiteness of primitive points breaking down happens for every curve when the degree is large enough.

\begin{theorem}\label{thm:main}
    Let $X$ be a nice curve over a number field $K$. Suppose that $X$ has a divisor of degree $d$ with $d> 2g(X)$. Then $K(X)$ contains a function $f$ of degree $d$ such that $\overline{K}(f) \subseteq \overline{K}(X)$ is a primitive field extension. In particular, $X$ has infinitely many primitive points of degree $d$ over $K$.
\end{theorem}
Assume that $X(K) \neq \emptyset$, then $X$ has a divisor of every degree. In particular, there is the following corollary.
\begin{corollary}
If $X(K) \neq \emptyset$ and $d > 2g(X)$, then $X$ has infinitely many primitive points of degree $d$ over $K$.
\end{corollary}

In a recent article by Neftin an Zieve the classification of primitive functions on curves of large enough degree has been completed. To be more precise they prove the following, where we omitted the results for $g=0$ and $g=1$ for simplicity.
\begin{theorem}[{\cite{neftin2024monodromygroups}[Thm. 1.1]}]\label{neftinzieve}
Let $g \geq 2$ be an integer.  Then there exist a constant $N_g$ such that for all curves $X$ over $\C$ and all $f \in \C(X)$ of degree $d \geq N_g$ the following holds: If $\C(f) \subseteq \C(X)$ is a primitive field extension then the Galois group of $\Gal(\C(X)/\C(f))$ is either $A_d$ or $S_d$.
\end{theorem}

For a function $f \in K(X)$ one has $\Gal(\overline K(X)/\overline K(f)) \cong \Gal(\C(X)/\C(f))$. In particular, $\overline K(f) \subset \overline K(X)$ is a primitive field extension if and only if  $\C(f) \subset \C(X)$ is.
Combining \Cref{thm:main} with \Cref{neftinzieve} gives the following corollary:

\begin{corollary}\label{cor:main}
	Let $X$ be a nice curve over a number field $K$ of genus $g \geq 2$. Suppose that $X$ has a divisor of degree $d$ with $d \geq \max(2g(X)+1, N_g)$; where $N_g$ is as in \Cref{neftinzieve}. Then $K(X)$ contains a function $f$ of degree $d$ such that $\Gal(\overline K (X)/\overline K(f))$ is either $A_d$ or $S_d$. In particular, $X$ has infinitely many points of degree $d$ with Galois group $A_d$ or $S_d$.
\end{corollary}

\Cref{sec:background} contains some basic definitions and facts regarding primitive points. \Cref{sec:contracting_divisors} contains some results useful for studying how primitive functions behave inside a fixed Riemann-Roch space. Finally, \cref{sec:proof_main_thm} is dedicated to the proof of \Cref{thm:main} which relies heavily uses the tools developed in \Cref{sec:contracting_divisors}. In fact, \Cref{thm:main} is deduced from a slightly stonger statement, see \Cref{prop:main}.

\subsection{Acknowledgements}

I would like to thank Maleeha Khawaja for introducing me to the subject of primitive points on curves and suggesting improvements to an earlier version of this paper.
The author was supported by the Croatian Science Foundation under the project no. IP-2022-10-5008 during this research.

\section{Background theory}\label{sec:background}

The goal of this section is to quickly recall some of the background theory on primitive points as developed by Khawaja and Siksek in \cite{KhawajaSiksek,khawaja2024singlesource}.

\begin{definition}
Let $K$ be a field, a finite field extension $L$ of $K$ is called primitive if $K$ and $L$ are the only sub fields of $L$ containing $K$. If $L$ is not a primitive extension of $K$ then it called imprimitive. 
\end{definition}

\begin{definition}
Let $X$ be a nice curve over a field $K$. Let $P$ be a point on $X$. The point $P$ is called a primitive point if $K(P)$ is primitive extension of $K$ and $P$ is called imprimitive if $K(P)$ is imprimitive.
\end{definition}

If $K$ is a perfect field and $L$ is a finite extension of $K$. Then $\widetilde L$ denotes the normal closure of $L$ over $K$. We will use $\Gal(L/K) := \Gal(\widetilde L/K)$ as shorthand notation of the Galois group of the normal closure of $L$.  

\begin{lemma}\label{lem:primitivity_galois}Let $K$ be a perfect field and $L$ be a finite field extension of $K$. Let $\tilde L$ be the normal closure of $L/K$. Then $L$ is primitive over $K$ if and only if the action of $\Gal(L/K) := \Gal(\widetilde L/K)$ on $Hom_K(L, \widetilde L)$ is primitive. 
\end{lemma}
\begin{proof}
The proof is essentially the same as that of \cite{KhawajaSiksek}[Lem. 28]. Although one has to use the fact that if one writes $L = K(\theta)$ and lets $\theta_1, \ldots \theta_d$ be the conjugates of $\theta$ in $ \tilde L$ then $\set{\theta_1, \ldots \theta_d}$ and $Hom_K(L, \widetilde L)$ are isomorphic as sets with a $\Gal(L/K)$ action. Where $f \in Hom_K(L, \tilde L)$ corresponds to $f(\theta) \in \set{\theta_1, \ldots \theta_d}$ under this isomorphism. 
\end{proof}

\begin{lemma}\label{lem:primitive_function_to_points}
    Let $X$ be nice curve over a number field $K$, suppose that $f \in K(X)$ is a function of degree $d$ such that $K(f) \subseteq K(X)$ is a primitive extension, then $X$ has infinitely many primitive points of degree $d$ over $K$.
\end{lemma}

\begin{proof}
Let $\widetilde{K(X)}$ be the normal closure of $K(X)$ over $K(f)$ and let $$G := \Gal(K(X)/K(f)) = \Gal(\widetilde{K(X)}/K(f)).$$ Then by \cite{Serre:Lectures-on-MW}[\S 9.2 Prop. 2] there is a thin set $S \subset \P^1(K)$ such that $$\Gal(K(f^{-1}(t))/K)=G$$ for all $t \in \P^1(K) \setminus S$, in particular if $t \in \P^1(K) \setminus S$ and $P \in f^{-1}(t)$ then $Hom_K(K(P), \widetilde{K(P)}) = Hom_K(f)(K(X), \widetilde{K(X)})$ as sets with a $\Gal(K(P)/K) = \Gal(K(f^{-1}(t))/K)=G$ action. In particular, \cref{lem:primitivity_galois} implies that for $t \in \P^1(K) \setminus S$ and $P \in f^{-1}(t)$ one has that $K(P)/K$ is primitive if and only if $K(X)/K(f)$ is primitive. The lemma now follows from Hilbert's irreducibly theorem \cite{Serre:Lectures-on-MW}[\S 9.6 Hilbert's Thm] which states that $\P^1(K) \setminus S$ is infinite.
\end{proof}

\section{Contracting divisors and imprimitive functions} \label{sec:contracting_divisors}
\begin{lemma}\label{lem:nice_divisor}
    Let $X$ be a nice curve over a number field $K$. Suppose that $X$ has a divisor $D_0$ of degree $d$ with $d> 2g(X)-2$. Then $X$ has an effective divisor $D$ of degree $d$ where every point in the support of $D$ occurs with multiplicity 1 in $D$ and such that furthermore  $H^0(X,\O(D))$ contains a function of degree exactly $d$.
\end{lemma}
\begin{proof}
Let $D_0$ be a divisor of degree $d$ on $X$. By the Rieman-Roch Theorem and $d \geq g(X)$ one has that $D_0$ is linearly equivalent to an effective divisor $D_1$.  By Riemann-Roch and $d > 2g(X)-2$ one has $\dim H^0(X,\O(D_1)) = d - g(X) +1 >0$ and $\dim H^0(X,\O(D_1 - D')) < d - g(X) +1$ for any nontrivial effective divisor $D'$, so that $H^0(X,\O(D_1))$ has to contain a function $f$ of degree exactly $d$.
The effective divisor $D$ can be obtained as $f^*(x)$ where $x \in \P^1(K)$ is a point above which $f$ is not ramified. The sought function in $H^0(X,\O(D))$ can be found by composing $f$ with an automorphism of $\P^1(K)$ that moves $x$ to $\infty$.
\end{proof}

\begin{definition}
Let $f: X \to C$ be a morphism between nice curves over a field $K$, and $D$ a divisor of degree $d$ on $X$. Then $f$ is called a contracting morphism for $D$ if 
\begin{enumerate}
    \item $1 < \deg f < d$ and 
    \item there exists a divisor $D'$ of degree $d/(\deg f)$ on $C$ such that $D = f^{*}(D')$.
\end{enumerate}
We define $\Contr(D)$ to be the set of isomorphism classes of contracting morphisms for $D$ defined over $K$. For a field extension $L$ of $K$ we use  $\Contr(D_L)$ to denote those that are defined over $L$. For an integer $g$ we let $\Contr_{g}(D_L) \subset \Contr(D_L)$ be the subset where the genus of the target curve $C$ is restricted to be equal to $g$.
\end{definition}

Note that we call two contracting morphisms $f: X \to C$ and $f' : X \to C'$ isomorphic if there is an isomorphism $g: C \to C'$ such that $f' = g \circ f$.

The main motivation for the above definition is the following lemma, which allows us to get a better understanding of the imprimitive functions in $H^0(X,\O(D))$.

\begin{lemma} \label{lem:imprimitive_factors}
    Let $X$ be a nice curve over a field $K$ and $D$ a divisor on $X$ and let $f$ be a function on $X$ whose pole divisor is $D$. If $K(f) \subset K(X)$ is a imprimitive extension then $f$ factors via a morphism in $\Contr(D)$.
\end{lemma}
\begin{proof}
If $K(f) \subset K(X)$ is not primitive then there is a field $L \subsetneq K(X)$ such that $K(f) \subsetneq L$. This field $L$ is the function field of some curve $C$ and this gives a factorisation of $f$ as $X \stackrel{g}{\to} C \stackrel{h}{\to} \P^1_K$ with $1 <\deg g < d$. Since $D = f^*(\infty)$ by assumption we can take $D':=h^*(\infty)$ to get $D = f^*(\infty) =g^*(h^*(\infty))=g^*(D')$. This shows that $f$ factors via $g$ with $g \in \Contr(D)$.
\end{proof}

\begin{proposition}\label{prop:contracting_iso_unique}
Let $f: X \to C$ and $f' : X \to C'$ be two contracting morphisms for some divisor $D$ on $X$, if $f$ and $f'$ are isomorphic then then there exists a unique $g: C \to C'$ such that $f' = g \circ f$. In particular, if $f$ and $f'$ become isomorphic over $K^{sep}$ they were already isomorphic over $K$.
\end{proposition}

\begin{proof}
    $g$ is determined uniquely by the map it induces between $C(\overline K) \to C'(\overline K)$. The morphism $f$ induces a surjective map $X(\overline K) \to C(\overline K)$. In particular, the relation  $f' = g \circ f$ uniquely determines the value $g$ should take on every $C(\overline K)$ point. 
    Finally, isomorphisms over $K^{sep}$ descent to isomorphisms over $K$ because they are unique and hence Galois invariant.
\end{proof}

\begin{lemma}\label{lem:contr_finite}
Let $X$ be a nice curve over a perfect field $K$ and $D$ a divisor on $X$ such that all points in the support of $D$ occur with multiplicity $1$ in $D$. Then $\Contr(D)$ is finite.
\end{lemma}
\begin{proof}
Since $K$ is perfect one has $K^{sep}=\overline K$. By \cref{prop:contracting_iso_unique} we have that $\Contr(D)$ injects in to $\Contr(D_{\overline K})$, so we can assume that $K$ is algebraically closed.
Since $X$ cannot map to curves of higher genus than $g(X)$, one has $$\Contr(D) = \bigcup_{g=0}^{g(X)-1} \Contr_{g}(D).$$

For $g \geq 2$ we have that $X$ only admits finitely many morphisms to curves of genus $g$ by the de Franchis theorem, see \cite{Kani1986}[Thm 3] for an explicit bound on the number of such morphisms. So in this case $\Contr_{g}(D)$ is finite. 

Now to deal with the case $g=1$. It is know that $X$ admits only finitely many morphisms of degree $< d$ to genus 1 curves up to isomorphisms on the target, see \cite{Kani1986}[Thm 4 and the corollary after it] for explicit bounds on the number of such morphisms. So $\Contr_{1}(D)$ is finite as well. 

It remains to show that $\Contr_{0}(D)$ is finite. Note that since we assumed that $K$ is algebraically closed we have that any genus $0$ curve is isomorphic to $\P^1_K$. Let $f \in \Contr_{0}(D)$ and $D'$ a divisor on $\P^1_K$ such that $D = f^{*}(D')$. Because $\deg f < d$ we have that $\deg D' = d / \deg f > 1$. Because $D = f^{*}(D')$ and all points in $D$ occur with multiplicity $1$, the same holds for $D'$. In particular, $\deg D' >1$ implies that there are at least two points $c_0$ and $c_\infty$ in the support of $D'$. By composing $f$ with an automorphism of $\P^1_K$ we can ensure that $c_0 = 0, c_\infty=\infty$. After applying this automorphism $f$ induces a linear equivalence between $D_0 := f^{*}(0)$ and  $D_\infty := f^{*}(\infty)$. Furthermore $D_0$ and $D_\infty$ are effective divisors that are smaller then $D$. The isomorphism class of the function $f$ is determined by $D_0$ and  $D_\infty$. The above argument shows that $\#\Contr_{0}(D)$ is at most the number of pairs of effective divisors that are smaller then $D$. Since the set of effective divisors smaller then $D$ is finite, so is $\Contr_{0}(D)$.
    
\end{proof}

\begin{definition}
For a divisor $D$ of degree $d$ on $X$ let $\P(D)$ be the fiber above $\O(D)$ under $X^{(d)} \to \Pic^d X$. 
\end{definition}

Note that $\P(D)$ parametrises all effective divisors of degree $d$ that are linearly equivalent to $D$, and is isomorphic to $\P^n_K$ with $n :=\dim H^0(X, \O(D))-1$. In fact, for any field extension $L$ of $K$ we have $\P(D)(L) = ( H^0(X_L, \O(D)) \setminus 0)/L^*$. In other texts $\P(D)(L)$ is sometimes denoted by $|D_L|$.

\begin{remark}\label{rem:points_as_functions}
Note that  $H^0(X_L, \O(D))$ are the functions on $X_L$ whose pole divisor is at most $D$. So the identification $\P(D)(L) = ( H^0(X_L, \O(D)) \setminus 0)/L^*$ allows one to see points on $\P(D)$ as functions up to scalar multiplication. In particular if a property of functions is invariant under scalar multiplication it makes sense to look at the set of functions $\P(D)(L)$ in with that property.
\end{remark}

\begin{proposition}\label{prop:dimension_comparison}
    Let $f: X \to C$ be a morphism of degree $e > 1$, $D'$ an effective divisor on $C$ and $D := f^*(D')$. Assume that $\deg D > 2g(X)$, then $\dim \P(D) > \dim \P(D')$.
\end{proposition}

\begin{proof}
    Let $d := \deg D$ and $d' := \deg D'$ so that $d=ed'$. Since $d > 2g(X)-2$ we have $\dim \P(D) = \dim H^0(X, \O(D)) = d-g(X)$ by Riemann-Roch. Now $2g(X)-2 \geq e(2g(C)-2)$ by Riemann-Hurwitz. In particular $d' = d/e > 2g(X)/e \geq 2g(C)-2$ so $\dim \P(D') = d'-g(C)$.
    The proposition now follows from $\dim \P(D) -\dim \P(D') = d-g(X) - d' + g(C) \geq d(1-1/e) -g(X)  >0$
\end{proof}

\section{Proof of \Cref{thm:main}} \label{sec:proof_main_thm}

By \cref{lem:nice_divisor} we can find a divisor $D$ of degree $d$ on $X$, where all points in the support of $D$ occur in $D$ with multiplicity 1.
Our goal is to show that $H^0(X, \O(D))$ contains a function of degree $d$ such that $\overline K(f) \subset \overline K(X)$ is a primitive field extension.  This easily follows from the following proposition. With this function $f$ in hand, the existence of infinitely many primitive points of degree $d$ follows from \cref{lem:primitive_function_to_points}.

\begin{proposition}\label{prop:main}
	Let $X$ be a nice curve over a number field $K$. Suppose that $X$ has an effective divisor $D$ of degree $d>2g(X)$ such that every point in the support of $D$ occurs with multiplicity 1 in $D$ . Then there is a dense open subset $U \subset \PP(D)$ such that all functions $f \in U(K)$ are of degree $d$ and $\overline K(f) \subseteq \overline K(X)$ is primitive and of degree $d$.
\end{proposition}

\begin{proof}
Let $S \subset \P(D)$ be the locus where the degree of the function is smaller then $d$. Then $S$ is a closed variety, and by Riemann-Roch and $d > 2g(X)$ we know that $S \neq \P(D)$. Let $f \in S(\overline K) \subseteq \P(D)(\overline K)$ be such that $\overline K(f) \subseteq \overline K(X)$ is an imprimitive extension. Then  by \Cref{lem:imprimitive_factors} $f$ factors via function in $\Contr(D_{\overline K})$. For a $g \in \Contr(D_{\overline K})$ let $D_g$ denote the effective divisor such that $g^*(D_g) = D$. Then $f$ is contained in $$T := \bigcup_{g \in \Contr(D_{\overline K})} g^*(\P(D_g))\subset \P(D).$$ By \Cref{lem:contr_finite} the set $\Contr(D)$ is finite, showing that $T$ is a finite union of a closed sub varieties of $\P(D)$. By \Cref{prop:dimension_comparison} all the varieties making up $T$ have dimension less then that of $\P(D)$. Note that $T$ is stable under the action of $\Gal(\overline K/K)$ which allows one to see $T$ as a (possibly reducible) variety over $K$.

So $U := \P(D) \setminus (S \cup T)$ is the open and dense subvariety that is asserted to exsist by the proposiition.

\end{proof}

\bibliographystyle{alpha}
\bibliography{references.bib}{}

\end{document}